\documentclass[12pt,a4paper]{article}

\usepackage[left=2.54cm, right=2.54cm, top=2.54cm, bottom=2.54cm]{geometry}

\usepackage{needspace}

\usepackage{tikz}
\usepackage{rotating}
\usepackage{caption}

\usepackage{setspace}

\usepackage{enumerate}

\usepackage{graphicx}
\usepackage{amsfonts,amssymb,eucal}
\usepackage{booktabs}
\usepackage{bm}
\usepackage[reqno]{amsmath}
\usepackage[final,tracking=true,kerning=true,spacing=true,factor=1100,stretch=10,shrink=10]{microtype}
\usepackage{amsthm}

\newtheorem{Proposition}{Proposition}
\newtheorem{TT}{Theorem}

\providecommand{\keywords}[1]{\small{\textit{Keywords:}} #1}

\usepackage{longtable}
\usepackage{siunitx}
\makeatletter
\g@addto@macro\appendix{
    \renewcommand{\thetable}{\Alph{section}.\arabic{table}}
}
\makeatother

\title{On the Vertex Growth Rate \\
of Integer Partition Polytopes \\
and Their Non-Vertex Structure}

\author{Vladimir A. Shlyk \\ \texttt{v.shlyk@gmail.com}}

\begin{document}

\onehalfspacing

\date{\vspace{-5ex}}
\maketitle

\begin{abstract}
The Hardy-Ramanujan formula states that the asymptotic growth of the number of integer partitions of n is 
$\exp(\Theta(\sqrt{n}\,))$. We prove that the number of vertices of the integer partition polytope is bounded by $\exp\left(\mathcal{O}(\ln^2 n)\right)$. This implies that all partitions of $n$ can be generated
from a negligible subset of vertex partitions using convex combinations.

We also present new computational data on the structure of non-vertex partitions,
which cause the huge gap between the number of vertices and the total number of partitions.
\end{abstract}

\section{Introduction}

\subsection{Partition Polytope, Its Vertices, and the Main Result}

The idea of investigating the set of partitions of a number $n \in \mathbb{Z}_+$ from a polyhedral perspective
was proposed in \cite{S96} and further developed in \cite{S05}. It is based on the construction of the partition polytope
$P_n \subset \mathbb{R}^n$, defined as the convex hull of all partitions of $n$.
This polytope is an analogue of the corner polyhedron introduced by Gomory \cite{G},
which plays a key role in the mathematics of integer programming.

The polytope $P_n$ is defined as the convex hull of all partitions of $n,$ where partitions are identified with their multiplicity vectors $ x=(x_1, x_2, \ldots, x_n) $:
$$
P_n := \mathrm{conv\,} \{ x = (x_1, x_2, \ldots, x_n) : 1x_1 + 2x_2 + \dotsb + nx_n = n \}.
$$

In \cite{S05} we characterized the facets of $P_n$ and then focused on its vertices.
Unlike most polytopes associated with ``natural`` combinatorial objects (such as permutations, graphs, or matroids)
--- where all the underlying objects are vertices --- this turned out not to be the case for the polytope $P_n$.
The first time when we meet a partition, that is not a vertex, is already for $n=4.$
The multiplicity vector $x=(2, 1, 0, 0)$ of the partition $4=1+1+2$ is not a vertex of $P_n$
because it is the half-sum of two other partitions $x=\frac12 [(0, 2, 0, 0)+(4, 0, 0, 0)].$
Consequently, for any $n,$ the vertices of $P_n$ form a kind of basis for the set of all partitions
under convex combinations, just as permutation matrices form a basis for the Birkhoff polytope.
So, the vertices have become a subject of particular interest.

We computed the vertices of $P_n$ up to $n = 105$ \cite{VS}, see sequence A203898 of their numbers $v(n)$
in the Online Encyclopedia of Integer Sequences (OEIS) at [https://oeis.org/A203898].
The computation data indicated that actually far from all partitions of $n$ are vertices;
the graph of $v(n)$ lies significantly lower than the graph of the function $p(n),$ the total number of partitions of~$n.$
For instance, at $n = 105$, the total number of partitions reaches $p(105) = 342\,325\,709$,
whereas the vertex basis consists of only $v(105) = 94\,520$ vertices --- less than $0.03\%$ of $p(105)$.
In \cite{S23} we conjectured that the growth rate of $v(n)$ is lower than that of $p(n),$
which --- as is known from the famous Hardy-Ramanujan formula --- is $\exp(\Theta(\sqrt n\, ))$ \cite{HR},
independently proved by Uspenskij \cite{U}.
In this paper, we show that this hypothesis is correct. The main result of the paper is
\begin{TT} \label{Bound}
The asymptotic bound holds:
\[
v(n) = \exp\left(\mathcal{O}(\ln^2 n)\right).
\]
\end{TT}

The path to the proof is indirect. Although in the joint work with 
Onn \cite{OS},
we proved that the problem of recognizing vertices of $P_n$ is solvable in polynomial time, that was done
with the use of linear programming methods, since no combinatorial characterization of vertices is known yet.

However, we have several necessary and several sufficient conditions for the vertices and, specifically,
a combinatorial characterization of partitions that are convex combinations
of two other partitions of the same number \cite{S14}.
Some consequences of this criterion served as the basis for the proof of Theorem \ref{Bound}.

\subsection{Computational Results on the Structure of the Set of Non-Vertex Partitions of $\bm{n}$}

In \cite{S23} we classified the non-vertex partitions of $n$ in terms of the minimum number $\xi$ of other partitions of which they are a convex combination. We call this number the convexity index of the partition. As a result of the vertex computation, we found that the majority of non-vertex partitions are of index 2. We also discovered numerous partitions of index 3 and a couple of index 4. But it left unknown whether partitions of index 5 and higher exist at all. Nevertheless, we conjectured that for any $\xi \geqslant 2$, a sufficiently large $n$ has partitions of index $\xi$~\cite{S23}.

We performed new calculations to verify the existence of index-5 partitions.
We found a significant number of such partitions for certain values of $n \geqslant 80$.
This partially confirmed the validity of the hypothetical structure of the set of partitions
of a given number proposed in \cite{S23}. In addition to the partitions themselves,
we recorded the coefficients of their convex representations in terms of the vertices
of one-point clean lattice simplices.
This establishes a bridge between the combinatorial geometry of partition polytopes
and the classification of toric singularities~\cite{K}.

The structure of the paper is as follows:
in Section 2, we provide a theoretical justification for the quasipolynomial upper bound
$\exp\left(\mathcal{O}(\ln^2 n)\right)$ on the number of vertices $v(n)$.
In Section~3, taking into account the significant gap between $p(n)$ and $v(n)$,
we describe the structure of the set of partitions of $n$ in greater detail
and present new computational data on non-vertex partitions
(specifically, on finding partitions of convexity index $\xi = 5$).

We use the following notation throughout the paper.
$x \vdash n$ means that $x$ is a partition of the number $n$.
We denote by $\mathcal{P}(n)$ the set of partitions of $n$, and by $\operatorname{Vert}(P_n)$
the set of vertices of the polytope $P_n$.
For $x \vdash n$, we call the set $S(x) := \{ i \in \{1, \dots, n\} : x_i > 0 \}$ the support of $x$,
and $d := |S(x)|$ is the size of the support.

\section{Quasipolynomial Growth Rate of $\bm{v(n)}$: \\ Proof of Theorem \ref{Bound}}

\subsection{Na\"\i ve Bound}

We can try to obtain an asymptotic upper bound on $v(n)$ using the inequality
$$v(n) \leqslant p(n) - p(n-4),$$
which holds because every partition $y$ of $n-4$, augmented with two parts $1$ and one part $2$,
is a convex combination of two partitions. However, calculating the right-hand side of this inequality
merely establishes that the relative density of the number of vertices tends to zero:
$$\lim_{n\to\infty} \frac{v(n)}{p(n)} = 0.$$

This demonstrates that the difference between the graphs of $v(n)$ and $p(n)$
is not merely a visual artifact, but a structural property: vertices constitute an asymptotically negligible fraction
of the set of all integer partitions of $n$.
Nevertheless, this elementary calculation preserves the exponential growth order of $\exp(\mathcal{O}(\sqrt{n}\,))$ for $v(n)$,
whereas the true divergence between $v(n)$ and $p(n)$ is far more substantial.
To obtain a significantly sharper bound, we must exploit the fine combinatorial structure of the vertices of $P_n$.
As established below, this leads to a quasipolynomial upper bound.

\subsection{Proof of the Quasipolynomial Bound on $\bm{v(n)}$}

\subsubsection{Key Inequalities}

The following necessary condition for $x\vdash n$ to be a vertex of $P_n$ was proved in \cite{S14}, Theorem 7:
\begin{equation}\label{Product}
\prod_{i \in S(x)} (x_i + 1) \leqslant n + 1.
\end{equation}
As a consequence, we obtain an exact upper bound on the number of distinct part sizes in any vertex:
\begin{equation}\label{NumParts}
d \leqslant \lfloor \log_2(n + 1)\rfloor.
\end{equation}
To make the presentation self-contained, we recall the proof of these inequalities. Let $x \in \operatorname{Vert}(P_n)$.
Consider all submultisets of parts of the partition $x$. In coordinate notation,
such a submultiset is defined by a vector $u = (u_1, \ldots, u_n)$, where $0 \leqslant u_i \leqslant x_i$. Its weighted sum is given by
$$\sigma(u) = 1 u_1 + 2 u_2 + \dotsb + n u_n .$$
For each index $i$, there are exactly $x_i + 1$ choices for $u_i$ (from $0$ to $x_i$).
Therefore, the total number of such submultisets (vectors $u$) equals
$$\prod_{i \in S(x)} (x_i + 1).$$
The next assertion is a half-statement of the criterion for a partition to be a convex combination of two other partitions:
\begin{center}
for $x \in \operatorname{Vert}(P_n)$, all numbers $\sigma(u)$ are distinct.
\end{center}
Indeed, suppose to the contrary that there exist two distinct submultisets $u \neq w$ such that $\sigma(u) = \sigma(w)$.
Define coordinatewise
$$a = (u - w)_{+}, \qquad b = (w - u)_{+},$$
where $(\cdot)_+$ denotes the positive part. Then $a$ and $b$ have disjoint supports, $a \leqslant x$, $b \leqslant x$, and
\begin{equation}\label{Ravvo}
1 a_1 + 2 a_2 + \dotsb + n a_n = 1 b_1 + 2 b_2 + \dotsb + n b_n.
\end{equation}
Since $u \neq w$, the vectors $a$ and $b$ are distinct;
equality (\ref{Ravvo}) then ensures that both vectors contain nonzero components. Setting
$$y = x - a + b, \qquad z = x + a - b,$$
and noting that $a \leqslant x$ and $b \leqslant x$, we have $y, z \in \mathbb{Z}_{+}^n$. Furthermore, equality (\ref{Ravvo}) implies
$$1y_1 + 2y_2 + \dotsb + ny_n = 1z_1 + 2z_2 + \dotsb + nz_n = n,$$
so $y$ and $z$ are valid partitions of $n$. Moreover, $y$ and $z$ are distinct from $x$ and from each other, satisfying
$$x = \frac{y + z}{2}.$$
Thus, $x$ is the midpoint of the line segment connecting two distinct points in $P_n$, contradicting the hypothesis that $x$ is a vertex.

Consequently, for any vertex $x$, the mapping $u \mapsto \sigma(u)$ is injective.
Since each $\sigma(u)$ is an integer taking values from $0$ to $n$ inclusive, the codomain contains at most $n + 1$ elements.
This immediately establishes the foundational inequality (\ref{Product}).

Using (\ref{Product}), we can bound the number of distinct parts $d = |S(x)|$.
For each $i \in S(x)$, we have $x_i + 1 \geqslant 2$. Thus,
$$2^d \leqslant \prod_{i \in S(x)} (x_i+1) \leqslant n+1.$$
Taking the base-$2$ logarithm and observing that $d$ is an integer, we obtain $d \leqslant L := \lfloor \log_2(n + 1) \rfloor$, proving (\ref{NumParts}).

Inequality (\ref{NumParts}) serves as the primary technical tool for our subsequent enumeration of vertices.

\subsubsection{Upper Bound on $\bm{v(n)}$}

Every vertex $x$ of $P_n$ is uniquely determined by its support and its vector of multiplicities.
For a fixed number $d$ ($1 \leqslant d \leqslant L$) of distinct parts in $x$,
the support $S(x)$ is a $d$-element subset of $\{1, \ldots, n\}$. The total number of such subsets is
$$C(n,d) = \binom{n}{d}.$$
It remains to estimate the number of feasible multiplicity vectors for a given support.

Set $N = n+1$ and $L = \lfloor \log_2 N \rfloor$. Fix a support $S = \{s_1, \ldots, s_d\}$ with $s_1 < \dotsb < s_d$ and $d \leqslant L$, and let $m_1, \ldots, m_d$ denote the positive multiplicities of its elements. By (\ref{Product}), these satisfy
\begin{equation}\label{Eq4}
(m_1+1)(m_2+1) \dotsb (m_d+1) \leqslant N.
\end{equation}
Introducing variables $y_j = m_j+1 \geqslant 2$ and integers $r_j = \lfloor \log_2 (y_j) \rfloor \geqslant 1$ for $j=1,\ldots,d$, the relation $2^{r_j} \leqslant y_j$ together with (\ref{Eq4}) implies
$$2^{r_1+\dotsb +r_d} \leqslant y_1 y_2 \dotsb y_d \leqslant N,$$
which reduces to the linear constraint
\begin{equation}\label{Eq5}
r_1 + \dotsb + r_d \leqslant L.
\end{equation}

To bound the number of valid multiplicity vectors, we analyze the distribution of the tuple $(r_1, \ldots, r_d)$.
Since these integers represent the discrete logarithms of (multiplicities $+ 1$), they are positive, may coincide, and their relative order is significant.
For a fixed sum $s \geqslant d$, the number of positive integer $d$-tuples satisfying $r_1+\dotsb +r_d=s$ is given by $C(s-1,d-1)$.
Summing over all permissible values of $s$, the total number of such tuples satisfying (\ref{Eq5}) equals
\begin{equation}\label{Eq6}
\sum_{s=d}^{L} C(s-1,d-1) = C(L, d),
\end{equation}
by the standard hockey-stick identity for binomial coefficients.

Next, we count the number of multiplicity choices corresponding to a fixed tuple $(r_1, \ldots, r_d)$.
For each component $r_j$, the variable $y_j$ is bounded by
$$2^{r_j} \leqslant y_j < 2^{r_j+1},$$
yielding at most $2^{r_j}$ integer choices for $y_j$.
For the entire $d$-tuple, the total number of multiplicity configurations does not exceed
\begin{equation}\label{Eq7}
2^{r_1} \dotsb 2^{r_d} = 2^{r_1+\dotsb +r_d} \leqslant 2^{L} \leqslant N.
\end{equation}
Combining (\ref{Eq6}) and (\ref{Eq7}), we conclude that for any fixed $d$-element support, the number of multiplicity vectors satisfying condition (\ref{Product}) is bounded by
\begin{equation}\label{Eq8}
N C(L,d).
\end{equation}

Summing over all $C(n,d)$ possible $d$-element supports across all $1 \leqslant d \leqslant L$, we obtain the total bound from (\ref{Eq8}):
\begin{equation}\label{Eq9}
v(n) \leqslant N \sum_{d=1}^{L} C(n,d) C(L,d).
\end{equation}
Applying the symmetry relation $C(L,d) = C(L, L-d)$ and Vandermonde's identity,
$$\sum_{d=0}^{L} C(n,d) C(L,L-d) = C(n+L,L),$$
we refine (\ref{Eq9}) into an explicit combinatorial bound:
\begin{equation}\label{Eq10}
v(n) \leqslant (n+1) \bigl(C(n+L,L) - 1\bigr), \qquad L = \lfloor \log_2(n+1) \rfloor,
\end{equation}
where the subtraction of $1$ accounts for the empty-support term $d=0$ in the full Vandermonde sum.

We now evaluate the asymptotic growth rate of this bound using the standard estimate
$$C(M,L) \leqslant \left(\frac{\mathrm{e} M}{L}\right)^L.$$
Setting $M = n + L$, inequality (\ref{Eq10}) yields
\begin{equation}\label{Eq11}
\ln v(n) \leqslant \ln(n+1) + L\bigl(1 + \ln(n+L) - \ln L\bigr).
\end{equation}
Expanding the logarithmic terms,
$$\ln(n+1) = \ln n + \mathcal{O}(1/n), \qquad L = \frac{\ln n}{\ln 2} + \mathcal{O}(1),$$
$$\ln(n+L) = \ln n + o(1), \qquad \ln L = \ln \ln n - \ln \ln 2 + o(1),$$
we simplify (\ref{Eq11}) to
$$\ln v(n) \leqslant \frac{(\ln n)^2}{\ln 2} - \frac{\ln n \cdot \ln \ln n}{\ln 2} + \mathcal{O}(\ln n).$$
This establishes the quasipolynomial upper bound:
$$v(n) = \exp \bigl( \mathcal{O}(\ln^2 n) \bigr) = n^{\mathcal{O}(\ln n)},$$
completing the proof of Theorem \ref{Bound}.

\subsubsection{Comparison of $\bm{v(n)}$ with $\bm{p(n)}$}

It is insightful to compare the asymptotic growth rate of $v(n)$ with that of the total number of integer partitions $p(n)$.
By the Hardy--Ramanujan asymptotic formula,
$$p(n) \sim \frac{1}{4n\sqrt{3}} \exp\left(\pi \sqrt{\frac{2n}{3}}\,\right),$$
so that
$$\ln p(n) = \Theta(\sqrt{n}\,).$$
In contrast, our structural bound yields
$$\ln v(n) = \mathcal{O}(\ln^2 n),$$
a growth rate heavily dominated by the subexponential growth of $p(n)$.
Thus, vertices represent an asymptotically negligible portion of all partitions:
$$v(n) = o(p(n)) \quad \text{as } n \to \infty.$$

In fact, our analysis establishes a substantially stronger asymptotic separation.
Combining $\ln p(n) = \Theta(\sqrt{n})$ with $\ln v(n) = \mathcal{O}(\ln^2 n)$, we see that the logarithmic growth rates themselves are decoupled:
$$\frac{\ln v(n)}{\ln p(n)} \to 0 \quad \text{as } n \to \infty,$$
which can be expressed equivalently as
$$v(n) = p(n)^{o(1)}.$$
This stronger formulation underscores the fundamental contrast between the rapid exponential growth of integer partitions and the restricted quasipolynomial landscape of polytope vertices.

\section{Computational Results}
We have not obtained the exact asymptotics of $v(n)$ since a lower bound of the same order has not been obtained.
So, it has not been proven that the resulting quasipolynomial upper bound is close to the true one.
Theorem \ref{Bound}, proved in the previous section, shows that almost all partitions are not vertices.
After computing vertices of the polytope $P_n$ and examination of non-vertex partitions
we classified them using a new notion of the convexity index of a partition.
This is the minimum number $\xi, \xi \geqslant 2$, of partitions of which a given partition is a convex combination.
Let $C_{\xi}(n), \xi \geqslant 2$, be the set of partitions of n that are of the index $\xi$;
and let $C_\xi$ be the union of $C_\xi(n)$ over all $n.$
We have performed new computations that demonstrate where exactly the non-vertices "go",
how they are distributed among classes $C_2, C_3, C_4$, and $C_5$.

\subsection{Structure of the Set of Partitions of $\bm{n}$}

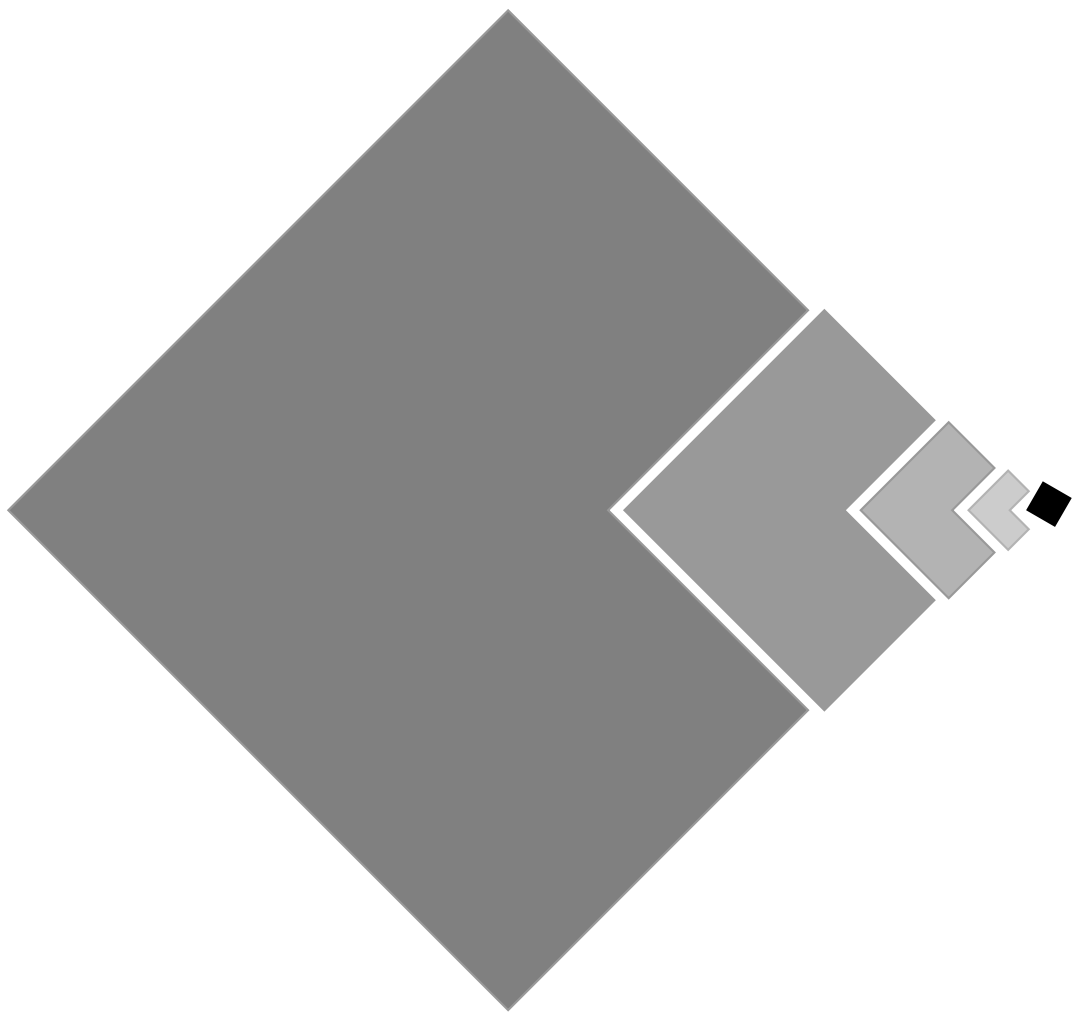
\begin{figure}[htbp]
\centering
\begin{tikzpicture}[scale=0.85]

    \def\S{5.5}
    \def\ppp{0.4}
    \def\pppp{0.44}
    \def\ppppp{0.45}
    \def\eps{0.09}
    \def\rotAngle{15}
    \pgfmathsetmacro\Sthree{\S*\ppp}
    \pgfmathsetmacro\Sfour{\Sthree*\pppp}
    \pgfmathsetmacro\Sfive{\Sfour*\ppppp}
    \def\kk{0.6}
    \def\kkk{0.55}
    \def\kkkk{0.52}
    \def\kkkkk{0.52}
    \pgfmathsetmacro\midA{\S * \kk}
    \pgfmathsetmacro\midThree{\Sthree * \kkk}
    \pgfmathsetmacro\midFour{\Sfour * \kkkk}
    \pgfmathsetmacro\midFive{\Sfive * \kkkkk}
    \begin{scope}[rotate=-45]
        \fill[black!50] (0,0) -- (0,\S) -- (\midA,\S) -- (\midA,\midA) -- (\S,\midA) -- (\S,0) -- cycle;
        \draw[thin, black!40] (0,0) -- (0,\S) -- (\midA,\S) -- (\midA,\midA) -- (\S,\midA) -- (\S,0) -- cycle;
        \pgfmathsetmacro\shiftCthree{\midA + \eps}
        \begin{scope}[shift={(\shiftCthree, \shiftCthree)}]
            \fill[black!40] (0,0) -- (0,\Sthree) -- (\midThree,\Sthree) -- (\midThree,\midThree) -- (\Sthree,\midThree) -- (\Sthree,0) -- cycle;
            \draw[thin, black!40] (0,0) -- (0,\Sthree) -- (\midThree,\Sthree) -- (\midThree,\midThree) -- (\Sthree,\midThree) -- (\Sthree,0) -- cycle;
            \pgfmathsetmacro\shiftCfour{\midThree + \eps}
            \begin{scope}[shift={(\shiftCfour, \shiftCfour)}]
                \fill[black!30] (0,0) -- (0,\Sfour) -- (\midFour,\Sfour) -- (\midFour,\midFour) -- (\Sfour,\midFour) -- (\Sfour,0) -- cycle;
                \draw[thin, black!40] (0,0) -- (0,\Sfour) -- (\midFour,\Sfour) -- (\midFour,\midFour) -- (\Sfour,\midFour) -- (\Sfour,0) -- cycle;
                \pgfmathsetmacro\shiftCfive{\midFour + \eps}
                \begin{scope}[shift={(\shiftCfive, \shiftCfive)}]
                    \fill[black!20] (0,0) -- (0,\Sfive) -- (\midFive,\Sfive) -- (\midFive,\midFive) -- (\Sfive,\midFive) -- (\Sfive,0) -- cycle;
                    \draw[thin, black!30] (0,0) -- (0,\Sfive) -- (\midFive,\Sfive) -- (\midFive,\midFive) -- (\Sfive,\midFive) -- (\Sfive,0) -- cycle;
                    \pgfmathsetmacro\shiftBlack{\midFive + \eps}
                    \begin{scope}[shift={(\shiftBlack, \shiftBlack)}, rotate=\rotAngle]
                        \fill[black] (0,0) rectangle (0.26, 0.26);
                    \end{scope}
                \end{scope}
            \end{scope}
        \end{scope}
    \end{scope}
\end{tikzpicture}
\caption{The hypothetical structure of $\mathcal{P}(n)$: the partition sets $C_2(n)$ -- $C_5(n)$ and the vertex set $\operatorname{Vert}P_n$
(the rotated black square on the right).}
\label{fig:partition_polytope_structure}
\end{figure}

Our previous computation established that the majority of non-vertex partitions have convexity index 2.
We also discovered numerous partitions of index 3 --- the first being $15=2\times 3+4+5$ ---
and two partitions of index 4: $36=7+8+9+12$ and $44=2\times 7+9+10+11$.
Our goal was to determine the number of vertices of the polytope.
Therefore, we did not focus on the specific characteristics of non-vertex partitions
and, in particular, did not aim to find partitions of higher indices.
After analyzing the behavior of the vertex-count function, we gained insight into the structure
of the set of partitions for each number, see \cite{S23}
and Figure \ref{fig:partition_polytope_structure}\footnote{After Anton Stankowski's
\textit{ Das Kleine ist sch\"on} (1956).}.

The entire square in Figure \ref{fig:partition_polytope_structure} corresponds to the set
$\mathcal{P}(n)$ of all partitions of $n$,
while the hexagons --- from left to right --- correspond to the sets of partitions of indices 2, 3, 4, and so on.
Finally, the tiny square on the right represents the set of vertices of the polytope $P_n$.
In \cite{S23}, we conjectured that for any index $\xi$ and sufficiently large $n$, partitions of $n$ of index $\xi$ exist.
Based on the computed data, we also conjectured that for each $n$ the sizes $|C_\xi(n)|$ decrease as $\xi$ increases.
So, the number of hexagons on the way to the last tiny square with vertices can be arbitrarily large.

The following theorem establishes a criterion for a partition to be of convexity index 2.
\begin{TT}[\cite{S14}, Theorem 5]\label{Kriterion}
A partition $x$ belongs to the class $C_2$ if and only if
there are two distinct submultisets of parts in $x$ of the same sum.
\end{TT}
However we know nothing about the classes $C_{\xi}$, $\xi >2$,
except two series of partitions in $C_3$
\cite{S14}, and that $C_4 \neq \emptyset$.

The following Proposition establishes relation between the convexity index of a partition and the size of its support.

\begin{Proposition}\label{d-geq-k}
Let $x$ be a partition of $n$ with support size $d = \vert{}S(x)\vert{}$. If $x \in C_k(n)$, then $d \geqslant k$.
\end{Proposition}

\begin{proof}
By definition, since $x \in C_k(n)$, it can be expressed as a convex combination
of a set of $k$ distinct partitions $Y = \{y^{(1)}, \dots, y^{(k)}\}$ of $n$:
\begin{equation}\label{eq:convex}
x = \sum_{j=1}^k \lambda_j y^{(j)}, \quad \sum_{j=1}^k \lambda_j = 1, \quad \lambda_j > 0,
\end{equation}
where $k$ is the minimum number of points among all such representations for $x$.

It follows from (\ref{eq:convex}) that a) for each $i \in S(x)$,
there exists at least one index $j \in \{1, \dots, k\}$ such that $y^{(j)}_i > 0$;
b) $y^{(j)}_i = 0$ for each $y^{(j)}$ and for every $i \notin S(x)$.
Thus, we have the exact equality
$$S(x) = \bigcup_{j=1}^k S(y^{(j)}).$$
Further on, viewing the partitions as points in \(\mathbb{R}^n\),
we identify them with their radius vectors in \(\mathbb{R}^d\) whose coordinates are indexed by \(S(x)\).

Let us show that the system of vectors $\{y^{(1)}, \dots, y^{(k)}\}$ is linearly independent.
Assume for contradiction that this system is linearly dependent.
By Carath\'eodory's theorem for cones (see \cite[Corollary 7.1i]{schrijver1986theory}),
any vector belonging to a cone generated by a set of vectors can be expressed
as a conic combination (i.e., a linear combination with non-negative coefficients)
of a linearly independent subset of those vectors.

Since all coefficients $\lambda_j$ in \eqref{eq:convex} are strictly positive, $x$ lies in the cone generated by $Y$.
The linear dependence of $Y$ implies that $x$ can be expressed as a convex combination of at most $k-1$ vectors in $Y$.
This contradicts the minimality of $k$. Thus, the vectors $y^{(1)}, \dots, y^{(k)}$ are linearly independent.

Since $\{y^{(1)}, \dots, y^{(k)}\}$
is a linearly independent system of $k$ vectors in the finite-dimensional space $\mathbb{R}^d$,
its cardinality cannot exceed the dimension of the space:
\[
k \leqslant \dim(\mathbb{R}^d) = d. \qedhere
\]
\end{proof}
Due to technical difficulties, the vertex computation had to be halted at $n = 105$.
We have now resumed the computations with a different goal:
to advance the verification of the hypothesis concerning the general existence of partitions of any convexity index, specifically by attempting to find at least one partition of index $5$.
Employing advanced algorithmic search techniques, we achieved results that surpassed our initial expectations.

\subsection{Detailed Computational Findings}

Before discussing the computational results, we should make an important clarification.
Initially, the computations covered all partitions for each $n$. However, starting at $n=76$,
we modified the primary algorithm to accelerate the process
and subsequently restricted the search to partitions with supports of size 5 and 6.
Recall that our main objective was to find at least one partition of index 5.
Hence, by Proposition \ref{d-geq-k}, we could expect the first such partition to have a support of size 5 or 6.

Consequently, the figures provided for partition counts up to and including $n=75$ reflect the complete picture,
as confirmed by verifying the equality
\begin{equation}\label{balance}
v(n) + |C_2(n)| + |C_3(n)| + |C_4(n)| = p(n)
\end{equation}
for all such $n$. However, for $n > 75$, it should be kept in mind
that the numbers of actually existing partitions of indices 3 and 4 are certainly greater than those presented,
and for partitions of index 5, they may also exceed the reported figures.

In total, we have found 38 partitions of index 5. Our first prize was the partition
$$80=9+16+17+18+20$$
of the number $n=80$. It turned out to be equal to one-fifth of the sum of the partitions
$80=4\times 20$, $80=5\times 16,$ $80=9+17+3\times 18,$ $80=9+3\times 17+20,$ and $80=3\times 9+17+2\times 18.$
In addition, it is a convex combination of partitions
$4\times 20,$ $5\times 16,$ $9+17+3\times 18,$ $9+3\times 17+20$, and $3\times 9+16+17+20$
with coefficients $1/6,$ $1/6,$ $1/3,$ $1/6$ and $1/6$ respectively.

However, for $n=100$ and $n=105$, we have already discovered nine such partitions for each.

Table \ref{tab : C5} presents the numbers of partitions of index $5$ on supports of size $5$ and $6$,
found for $80 \leqslant n \leqslant 105$.
Table \ref{tab:C2C3C4C5} provides complete computational data
for convexity indices $\xi = 2, 3, 4, 5$ for $n \leqslant 75$, along with available data for $n > 75$,
while Table \ref{tab : Partitions in C5} contains the explicit lists of these partitions of index~$5$.

\footnotesize
\begin{table}
$$\begin{array}{c|ccccccccccccc}
                        n & 80  & 81 & 82 & 83  & 84 & 85 & 86 & 87 & 88 & 89 & 90 & 91 & 92 \\ \hline
\vert{}C_5(n)\vert{}      &  1  &  0 &  0 &  0 &  0  &  2 &  1 &  0 &  0 &  0 &  1 &  2 &  1 \\ \hline
\hline
                        n & 93 & 94 & 95 & 96 & 97 & 98 & 99 & 100 & 101 & 102 & 103 & 104 & 105 \\ \hline
\vert{}C_5(n)\vert{}      &  0 &  1 &  4 &  0 &  0 &  0 &  0 &   9 &   1 &   2 &   2 &   2 &   9
\end{array}$$
\caption{Number of partitions of $n$ of convexity index 5.} \label {tab : C5}
\end{table}
\normalsize

From Table \ref{tab : C5} it is clear that with the exception of one of the first numbers in the list, $n=91,$
the maxima of this sequence occur at n multiples of 5.

\subsection{Statistics Across the Classes $\bm{C_3}$, $\bm{C_4}$, and $\bm{C_5}$}
We have registered plenty of partitions of index 3, and many partitions of index 4.
The statistics for these partitions are presented in Table \ref{tab:C2C3C4C5}.
It confirms the hypothesis proposed in \cite{S23} that these partitions are distributed across $n \in \mathbb{Z}_{+}$
according to the general rule: for each class $C_{\xi}$, they are particularly numerous when $n$ is a multiple of~$\xi$.
This could be seen especially clearly if graphs were plotted for these sequences.
This experimental observation also supports another hypothesis proposed in~\cite{S23}:
\begin{center}
\emph{the number of vertices of $P_n$ reaches local maxima at prime values of $n$.}
\end{center}

\subsection{Coefficients in Convex Representations}
That was far from all we accomplished. We recorded not only the partitions themselves
but also the coefficients in their clean representations as convex combinations.
By a clean representation of $x \in C_{\xi}$, we mean its convex representation
via the vertices of a $(\xi -1)$-dimensional one-point clean lattice simplex
in which $x$ is the unique strictly interior lattice point, and no other lattice points lie on its boundary.

It is known \cite{Re} that for a 2-dimensional lattice simplex (triangle) not associated with integer partitions,
all coefficients in the clean representation of the unique strictly interior lattice point are equal to $1/3$.
In the 3-dimensional case (which includes our class $C_4$), Reznick \cite{Re} established that
exactly seven distinct coefficient sets are theoretically possible.
Previously, only two partitions of index $4$ were known,
both yielding the symmetric coefficient set $\frac{(1, 1, 1, 1)}{4}$.
Our computations reveal that for index-$4$ partitions, in addition to this most frequent set,
one more theoretical possibility is realized, namely $\frac{(2, 1, 1, 1)}{5}$ (disregarding permutations).

Moving to a 4-dimensional clean simplex (with 5 vertices), Kasprzyk \cite{K} proved that
there exist $28\,686$ distinct coordinate types for the unique interior lattice point.
Because our simplices must reside entirely within the partition polytope $P_n$,
not all of these types are necessarily realized.
However, for the $5$-index partitions identified in our experiments,
we encountered three clean coefficient types (up to permutation):

$$t_5 = \frac{(1, 1, 1, 1, 1)}{5}, \quad t_6 = \frac{(2, 1, 1, 1, 1)}{6}, \quad t_{16} = \frac{(5, 4, 3, 3, 1)}{16},$$
all of which naturally belong to Kasprzyk's classification \cite{K}.
Specifically, types $t_5$ and $t_6$ appear regularly as clean representations,
whereas type $t_{16}$ was found in a single instance:
for the partition $100 = 7 + 17 + 20 + 25 + 31$,
which also exhibits representations of all the other coefficient types discussed below
($t_5, t_6, t_{10}, t_{12}, t_{15}, t_{21}$).

Furthermore, we identified convex representations via 4-simplex vertices
that contain one or more additional lattice points on the boundary (i.e., "non-clean" representations).
The corresponding coefficient types are:
$$t_{10} = \frac{(3, 2, 2, 2, 1)}{10}, \quad t_{12} = \frac{(4, 3, 2, 2, 1)}{12}, \quad t_{15} = \frac{(5, 3, 3, 3, 1)}{15}, \quad t_{21} = \frac{(7, 5, 4, 4, 1)}{21}.$$
Geometrically, the non-clean types $t_{10}$ and $t_{12}$ each correspond
to a single boundary point at the midpoint of an edge.
Type $t_{21}$ corresponds to a boundary point at the barycenter of a 3-face,
while type $t_{15}$ features two points on an edge combined with a 3-face barycenter.

Analytically, representations of types $t_{10}$, $t_{12}$, and $t_{15}$ can be reduced
to clean configurations by replacing a simplex vertex with an extra boundary point
(specifically, $t_{10}$ reduces to $t_5$, and $t_{12}$ reduces to $t_6$).
For type $t_{21}$, this reduction property remains unproven.

Now let us consider these coefficient types in detail for the 9 index-5 partitions found for $n = 105$.
These partitions are given by the following sets of distinct parts (each having multiplicity 1):
\begin{align*}
&p_1=\{10, 19, 20, 21, 35\}, & &p_2=\{13, 17, 19, 21, 35\}, & &p_3=\{15, 17, 21, 22, 30\}, \\
&p_4=\{16, 17, 21, 24, 27\}, & &p_5=\{17, 19, 20, 21, 28\}, & &p_6=\{1, 2, 9, 20, 24, 49\}, \\
&p_7=\{1, 2, 12, 21, 26, 43\}, & &p_8=\{1, 6, 9, 11, 28, 50\}, & &p_9=\{2, 5, 13, 17, 21, 47\}.
\end{align*}
In total, we found $25$ five-vertex representations for these partitions: $17$ clean and $8$ non-clean.
All 9 partitions except $p_5$ admit clean representations of both types $t_5$ and $t_6$,
whereas $p_5$ admits only type $t_6$. Specifically, four rich non-vertex partitions admit four representations
($t_5, t_6, t_{10}, t_{15}$); and one partition admits a unique clean representation of type $t_6$ (with no $t_5$).

\subsection{Computational Methods}
Our hunt for non-vertex partitions of index 5 was carried out in close collaboration with ChatGPT.

\noindent
\noindent\textbf{Algorithmic Approach.}

\noindent
In accordance with the refined search strategy, computations were restricted to supports of size 5 and 6.
The search exploited the structural property that a partition with support $S$
can only be a convex combination of partitions whose supports are subsets of $S$.

\noindent
For each integer $n$, the procedure unfolded as follows:
\begin{enumerate}
  \item \textbf{Support Filtering:} We enumerated subsets $S \subset \{1, 2, \dots, n\}$ of cardinality 5 and 6.
        Initial feasibility criteria, derived from the characterization of class $C_2$ (see Theorem~\ref{Kriterion}),
        were applied to filter out infeasible supports.
        This rapidly eliminated vast groups of supports during construction,
        significantly reducing the computational workload.
        For instance, for $n = 105$, out of $\binom{105}{5} + \binom{105}{6} = 1\,705\,904\,746$
        possible subsets, only $275\,948$ remained.
  \item \textbf{Partition Generation and Index Evaluation:} For each feasible support $S$,
        we enumerated all non-negative integer solutions to $s_1 y_1 + \dots + s_m y_m = n$, where $s_i \in S$.
        For $n = 105$, the number of these solutions per support ranged from $6$ to $1\,689$.
        Candidate partitions with full support $S$ were evaluated sequentially
        to determine their vertex status or convexity index.
  \item \textbf{Solver Integration:}
        \begin{itemize}
            \item To check whether a candidate is a vertex of $P_n$, we tested it against the convex hull of all partitions on $S$ using the \texttt{scipy.optimize.linprog} module.
            \item To determine the exact convexity index for partitions of index 4 or higher,
                we solved the minimization problem for the number of partitions in a convex combination using the \texttt{scipy.optimize.milp} module (where continuous variables represented decomposition coefficients and binary variables indicated the inclusion of specific partitions).
        \end{itemize}
  \item \textbf{Validity Verification:} When a partition was determined to have index 5, independent checks were performed to confirm the correctness of the computation.
\end{enumerate}

\noindent
\textbf{Data Consistency and Table Interpretation.}

\noindent
Regarding the data presented in Table~\ref{tab:C2C3C4C5} across different ranges of $n$:
\begin{itemize}
  \item \textbf{For $n \leqslant 75$,} the exhaustive counts strictly satisfy (\ref{balance}).
  \item \textbf{For $n \geqslant 76$,} while the entries reflect the restricted domain of supports 5 and 6 rather than exhaustive totals, their dynamical behavior and peaks appear to capture the general distribution pattern.
      This reinforces the hypothesis from \cite{S23} that the structure of these partitions
      is governed by the arithmetic factorization of $n$.
\end{itemize}

\section{Concluding Remarks and Open Problems}
The primary theoretical result of this paper establishes the quasipolynomial upper bound
$\exp\bigl(\mathcal{O}(\ln^2 n)\bigr)$
on the growth rate of the number of vertices of the integer partition polytope~$P_n$.
It follows that the size $v(n)$ of the vertex basis of the partition set under convex linear combinations
is negligibly small compared to the total number of partitions $p(n)$ (for instance, $v(105)/p(105) < 0.03\%$).
Moreover, the proof of Theorem \ref{Bound} implies that even the logarithm of $v(n)$ is small relative to $p(n)$.

The computational results provide new evidence supporting the perspective formulated in \cite{S23}
regarding the internal structure of the partition set.
We confirmed the existence of non-vertex partitions of convexity index 5 within the vast partition space
and obtained fresh evidence linking the counts of vertices and non-vertices of various indices
to the arithmetic structure of the integer being partitioned.
This leads to the conclusion that, from a polyhedral viewpoint, the structure of the set of integer partitions ---
objects purely additive in nature --- is governed by multiplicative properties, specifically by the factorization of $n$.

We cannot be sure that the obtained asymptotics for the number of vertices is exact,
as a lower bound of the same order is currently unknown.
Thus, determining the exact estimate for $v(n)$ remains an open problem.

The proof of Theorem \ref{Bound} relies on inequality (\ref{NumParts}),
which bounds the number of parts of distinct sizes in partitions
that are not convex combinations of two partitions.
Investigating unrelated problems, Ehrenborg and Readdy \cite{ER} arrived at partitions which,
following Theorem~\ref{Kriterion}, constitute the set $\mathcal{P}(n) \setminus C_2(n)$.
They termed them knapsack partitions; however, to our knowledge,
these objects have not yet been the subject of thorough study.
Geometrically (as can be seen from Figure \ref{fig:partition_polytope_structure}),
they correspond to the entire square minus the largest left hexagon.
Inequality (\ref{Product}) in essence holds for knapsack partitions.
Thus, although our primary goal concerned vertices, in proving our main result
we actually worked with the larger class of knapsack partitions, establishing a growth rate bound for them as well.

From the polyhedral standpoint, the structure of the partition set of $n$ forms a union of disjoint subsets:
$$\mathcal{P}(n) = C_2(n) \cup C_3(n) \cup C_4(n) \cup C_5(n) \cup \dots \cup \operatorname{Vert} P_n,$$
where the chain length $k$ of $\bigcup_{\xi \leqslant k} C_{\xi}(n)$ depends on $n$.
We conjecture that the distribution of partitions of $n$ across these subsets satisfies the following properties:
\begin{enumerate}[(i)]
    \item for each $\xi \geqslant 2$, the proportion of partitions in $C_{\xi}(n)$ is particularly large when $n$ is a multiple of $\xi$;
    \item the majority of partitions lie in $C_2(n)$, with the proportion being significantly higher for even $n$
        than for odd $n$, a phenomenon explained in \cite{S23}.
    \item for prime indices $\xi \geqslant 3$, the proportions falling into $C_{\xi}(n)$ exhibit a steady decline with $\xi$;
    \item the pattern for composite indices ($\xi = 4, 6, 8, 9, \dots$) requires further investigation.
          Nevertheless, we conjecture that a general monotonicity pattern holds with respect to divisibility:
          for instance, the proportion in $C_4(n)$ is smaller than that in $C_2(n)$,
          and the proportion in $C_6(n)$ is smaller than those in both $C_2(n)$ and $C_3(n)$.
          When comparing the proportions in $C_4(n)$ and $C_5(n)$,
          we expect the former to be larger than the latter.
          However, it is conceivable that for sufficiently large $n$, 
          the polytope~$P_n$ provides enough geometric room for 4-dimensional clean one-point simplices,
          potentially causing the proportion of partitions in $C_5(n)$ to exceed that in $C_4(n)$.
    \end{enumerate}

Notably, this distribution is consistent with --- and to some extent reinforces~---
the picture presented from a slightly different angle in~\cite{S23},
which examined how the number of vertices depends on the factorization of $n$.

Understanding the structure of the non-vertex set offers valuable insight for refining the upper bound
on the number of vertices, as well as for formulating new conjectures regarding its lower bound.
Another intriguing open problem is to determine the minimum value of $n$
for which the set of partitions with a specific convexity index $\xi$ is non-empty.
Currently, we know these values for $\xi = 2, 3, 4, 5$, which are $4, 15, 36, 80$, respectively.
Is it reasonable to expect that for $\xi = 6$,
the corresponding value of $n$ will turn out to be only slightly greater than $160$?

A natural direction for future research lies in the structural and computational analysis
of the generating subset of vertices.
It is known that all vertices of $P_n$ can be recursively constructed from a significantly smaller subcollection
of support vertices using two elementary combinatorial operations, see \cite{S14}.
For instance, at $n = 100$, where complete data on support vertices is available at [OEIS A203899],
the total number of partitions reaches $p(100) = 190\,569\,292$,
while the vertex basis consists of $v(100) = 59\,294$ vertices, and the number of support vertices is only $7\,228$
(representing less than $13\%$ of the vertices and a mere $0.0038\%$ of all partitions).
Obtaining an explicit structural characterization of support vertices
and establishing asymptotic bounds on their growth rate remain open problems of strong combinatorial interest.

There lies a wide horizon ahead for theoretical and experimental research
into the polyhedral structure of integer partition sets, potentially assisted by deeper, alternative methods.

\section{Acknowledgements}
The author is grateful to ChatGPT for collaborative algorithmic design, computational setup,
and assistance in shaping the proofs.

\appendix
\renewcommand{\thesection}{Appendix \Alph{section}}

\newpage

\section{Computational Data Tables}

\label{app:tables}

\begin{small}
\setlength{\tabcolsep}{10pt}
\begin{longtable}{
    r
    S[table-format=8.0, table-number-alignment = right]
    S[table-format=6.0, table-number-alignment = right]
    S[table-format=5.0, table-number-alignment = right]
    S[table-format=1.0]}
\caption{\spaceskip=0.7em {~}Complete computational data for convexity indices ${~~~~~~~~~~}\xi=2,\,3,\,4,\,5$ for $29 \leqslant n \leqslant 105$.\label{tab:C2C3C4C5}} \\
\toprule
$n$ & {$|C_2(n)|$} & {$|C_3(n)|$} & {$|C_4(n)|$} & {$|C_5(n)|$} \\
\midrule
\endfirsthead
\multicolumn{5}{c}{Table \thetable{} -- Continued from previous page} \\
\toprule
$n$ & {$|C_2(n)|$} & {$|C_3(n)|$} & {$|C_4(n)|$} & {$|C_5(n)|$} \\
\midrule
\endhead
\midrule
\multicolumn{5}{r}{{(Continued on next page)}} \\
\endfoot
\bottomrule
\endlastfoot
29 & 4131 & 1 & 0 & 0 \\
30 & 5256 & 6 & 0 & 0 \\
31 & 6295 & 3 & 0 & 0 \\
32 & 7884 & 1 & 0 & 0 \\
33 & 9479 & 18 & 0 & 0 \\
34 & 11722 & 1 & 0 & 0 \\
35 & 14047 & 11 & 0 & 0 \\
36 & 17296 & 10 & 1 & 0 \\
37 & 20623 & 6 & 0 & 0 \\
38 & 25142 & 4  & 0 & 0 \\
39 & 29942 & 29 & 0 & 0 \\
40 & 36299 & 12 & 0 & 0 \\
41 & 43081 & 11 & 0 & 0 \\
42 & 51950 & 31 & 0 & 0 \\
43 & 61439 & 17 & 0 & 0 \\
44 & 73668 & 10 & 1 & 0 \\
45 & 87040 & 62 & 0 & 0 \\
46 & 103748 & 16 & 0 & 0 \\
47 & 122149 & 31 & 1 & 0 \\
48 & 145155 & 59 & 1 & 0 \\
49 & 170487 & 54 & 1 & 0 \\
50 & 201710 & 27 & 1 & 0 \\
51 & 236380 & 119 & 0 & 0 \\
52 & 278616 & 29 & 3 & 0 \\
53 & 325705 & 58 & 2 & 0 \\
54 & 382735 & 80 & 1 & 0 \\
55 & 446379 & 108 & 1 & 0 \\
56 & 522753 & 81 & 1 & 0 \\
57 & 608499 & 213 & 0 & 0 \\
58 & 710464 & 67 & 2 & 0 \\
59 & 825116 & 115 & 1 & 0 \\
60 & 961126 & 185 & 8 & 0 \\
61 & 1113855 & 148 & 2 & 0 \\
62 & 1293807 & 96 & 3 & 0 \\
63 & 1496666 & 329 & 2 & 0 \\
64 & 1734274 & 136 & 10 & 0 \\
65 & 2002527 & 246 & 1 & 0 \\
66 & 2315287 & 314 & 4 & 0 \\
67 & 2668027 & 260 & 5 & 0 \\
68 & 3078114 & 160 & 26 & 0 \\
69 & 3541138 & 560 & 3 & 0 \\
70 & 4077157 & 279 & 7 & 0 \\
71 & 4681937 & 349 & 5 & 0 \\
72 & 5380381 & 507 & 13 & 0 \\
73 & 6168513 & 424 & 9 & 0 \\
74 & 7075312 & 283 & 11 & 0 \\
75 & 8098871 & 929 & 10 & 0 \\
76 & {--} & 248 & 38 & 0 \\
77 & {--} & 481 & 14 & 0 \\
78 & {--} & 419 & 17 & 0 \\
79 & {--} & 492 & 14 & 0 \\
80 & {--} & 341 & 30 & 1 \\
81 & {--} & 836 & 10 & 0 \\
82 & {--} & 418 & 22 & 0 \\
83 & {--} & 679 & 19 & 0 \\
84 & {--} & 816 & 47 & 0 \\
85 & {--} & 822 & 18 & 2 \\
86 & {--} & 568 & 33 & 1 \\
87 & {--} & 1383 & 29 & 0 \\
88 & {--} & 770 & 61 & 0 \\
89 & {--} & 1069 & 32 & 0 \\
90 & {--} & 1044 & 38 & 1 \\
91 & {--} & 1367 & 38 & 2 \\
92 & {--} & 936 & 100 & 1 \\
93 & {--} & 2176 & 34 & 0 \\
94 & {--} & 1084 & 48 & 1 \\
95 & {--} & 1789 & 48 & 5 \\
96 & {--} & 1968 & 109 & 0 \\
97 & {--} & 2000 & 59 & 0 \\
98 & {--} & 1538 & 68 & 0 \\
99 & {--} & 3190 & 46 & 0 \\
100 & {--} & 1751 & 171 & 9 \\
101 & {--} & 2617 & 87 & 1 \\
102 & {--} & 2732 & 91 & 2 \\
103 & {--} & 3018 & 85 & 2 \\
104 & {--} & 2380 & 166 & 2 \\
105 & {--} & 5088 & 102 & 9 \\

\end{longtable}

\begin{longtable}{r >{\raggedleft\arraybackslash}p{3.9cm} >{\kern-1.5cm}r}

\caption{Partitions in $C_5(n)$ with clean types, $80\leqslant n \leqslant 105$.} \label{tab : Partitions in C5} \\
\toprule
$n$ & Partition & Clean types \\
\midrule
\endfirsthead

\multicolumn{3}{c}{{\bfseries Table \thetable{} -- Continued from previous page}} \\
\toprule
$n$ & Partition & Clean types \\
\midrule
\endhead

\midrule
\multicolumn{3}{r}{{\itshape (Continued on next page)}} \\
\endfoot

\bottomrule
\endlastfoot

0 & $\{9, 16, 17, 18, 20\}$ & $t_5, t_6$ \\
85 & $\{9, 17, 18, 19, 22\}$ & $t_5, t_6$ \\
85 & $\{13, 15, 16, 17, 24\}$ & $t_5, t_6$ \\
86 & $\{1, 5, 7, 9, 23, 41\}$ & $t_5, t_6$ \\
90 & $\{2, 5, 10, 14, 18, 41\}$ & $t_6$ \\
91 & $\{13, 15, 18, 19, 26\}$ & $t_5, t_6$ \\
91 & $\{15, 16, 17, 19, 24\}$ & $t_5, t_6$ \\
92 & $\{15, 17, 18, 19, 23\}$ & $t_5, t_6$ \\
94 & $\{15, 17, 18, 19, 25\}$ & $t_5, t_6$ \\
95 & $\{13, 16, 19, 23, 24\}$ & $t_5, t_6$ \\
95 & $\{14, 17, 18, 19, 27\}$ & $t_5, t_6$ \\
95 & $\{15, 17, 18, 19, 26\}$ & $t_5$ \\
95 & $\{1, 5, 10, 17, 19, 43\}$ & $t_5, t_6$ \\
100 & $\{7, 17, 20, 25, 31\}$ & $t_5, t_6, t_{16}$ \\
100 & $\{9, 19, 20, 25, 27\}$ & $t_5, t_6$ \\
100 & $\{12, 20, 21, 22, 25\}$ & $t_5, t_6$ \\
100 & $\{14, 17, 20, 24, 25\}$ & $t_5, t_6$ \\
100 & $\{16, 18, 19, 20, 27\}$ & $t_5, t_6$ \\
100 & $\{1, 2, 11, 20, 25, 41\}$ & $t_5, t_6$ \\
100 & $\{1, 5, 7, 9, 30, 48\}$ & $t_5, t_6$ \\
100 & $\{1, 5, 8, 11, 27, 48\}$ & $t_5$ \\
100 & $\{1, 5, 10, 17, 19, 48\}$ & $t_5, t_6$ \\
101 & $\{15, 17, 20, 21, 28\}$ & $t_5, t_6$ \\
102 & $\{17, 18, 19, 21, 27\}$ & $t_6$ \\
102 & $\{2, 7, 10, 11, 25, 47\}$ & $t_5, t_5, t_6, t_6$ \\
103 & $\{17, 18, 19, 21, 28\}$ & $t_5$ \\
103 & $\{17, 18, 19, 22, 27\}$ & $t_5$ \\
104 & $\{9, 19, 23, 26, 27\}$ & $t_5, t_6$ \\
104 & $\{1, 5, 7, 9, 32, 50\}$ & $t_5, t_6$ \\
105 & $\{10, 19, 20, 21, 35\}$ & $t_5, t_6$ \\
105 & $\{13, 17, 19, 21, 35\}$ & $t_5, t_6$ \\
105 & $\{15, 17, 21, 22, 30\}$ & $t_5, t_6$ \\
105 & $\{16, 17, 21, 24, 27\}$ & $t_5, t_6$ \\
105 & $\{17, 19, 20, 21, 28\}$ & $t_6$ \\
105 & $\{1, 2, 9, 20, 24, 49\}$ & $t_5, t_6$ \\
105 & $\{1, 2, 12, 21, 26, 43\}$ & $t_5, t_6$ \\
105 & $\{1, 6, 9, 11, 28, 50\}$ & $t_5, t_6$ \\
105 & $\{2, 5, 13, 17, 21, 47\}$ & $t_5, t_6$ \\

\end{longtable}
\end{small}


\begin{thebibliography}{20}

\bibitem{ER}
Ehrenborg, R., Readdy, M. A. (2007). The M\"{o}bius function of partitions with restricted block sizes.
{\it Advances in Applied Mathematics} 39: 283--292. DOI: 10.1016/j.aam.2006.08.006.

\bibitem{G}
Gomory, R. E. (1969). Some polyhedra related to combinatorial problems.
{\it Linear Algebra and its Applications} 2: 451--558. DOI: 10.1016/0024-3795(69)90017-2.

\bibitem{HR}
Hardy, G. H., Ramanujan, S. (1918). Asymptotic formul{\ae} in combinatory analysis.
{\it Proceedings of the London Mathematical Society, Second Series} 17: 75--115.
DOI: 10.1112/plms/s2-17.1.75.

\bibitem{K}
Kasprzyk, A. M. (2013). Classifying terminal weighted projective space,
arXiv preprint arXiv:1304.3029 (Data accessible via the Graded Ring Database: http://www.grdb.co.uk).

\bibitem{OS}
Onn, S., Shlyk, V. A. (2015). Some efficiently solvable problems over integer partition polytopes.
{\it Discrete Applied Mathematics} 180: 135--140. DOI: 10.1016/j.dam.2014.08.015.

\bibitem{Re}
Reznik, B. (1986). Lattice point simplices.
{\it Discrete Mathematics} 60: 219--242. DOI: 10.1016/0012-365X(86)90015-4.

\bibitem{schrijver1986theory}
Schrijver, A.(1986).\emph{Theory of Linear and Integer Programming}, John Wiley \& Sons.

\bibitem{S96}
Shlyk, V. A. (1996). Polytopes of partitions of numbers.
{\it Proceecings of the National Academy of Sciences of Belarus. Phys.-Math. Ser.} 3: 89--92. (In Rusian).

\bibitem{S05}
Shlyk, V. A. (2005). Polytopes of partitions of numbers.
{\it European Journal of Combinatorics} 26: 1139--1153. DOI: 10.1016/j.ejc.2004.08.004.

\bibitem{DAN}
Shlyk, V. A. (2008). On the vertices of the polytopes of partitions of numbers.
{\it Doklady Nacionalnoj Akademii Nauk Belarusi} 52: 5--10. (In Russian).

\bibitem{S14}
Shlyk, V. A. (2014). Polyhedral approach to integer partitions.
{\it Journal of Combinatorial Mathematics and Combinatorial Computing} 89: 113--128.

\bibitem{S23}
Shlyk, V. A. (2023). Number of vertices of the polytope of integer partitions and factorization of the partitioned number. {\it Experimental Mathematics} 32(3): 457--466. DOI: 10.1080/10586458.2021.1980455.

\bibitem{VS}
Vroublevski, A. S., Shlyk, V. A. (2015). Computing vertices of integer partition polytopes.
{\it Informatics} 4: 34--48 (In Russian).

\bibitem{U}
Uspenskij, J. V. (1920). Les expressions asymntotiques des fonctions num\'eriques,
co\"{\i}ncidant dans les probl\`emes de dislocation des nombres en compos\'es.
{\it Bulletin de l'Acad\'emie des Sciences de Russie. VI s\'erie} 14: 199--218.

\end{thebibliography}
\end{document}